\documentclass[11pt]{article}
\usepackage{amssymb}
\usepackage{graphicx}
\usepackage{amsfonts}
\usepackage{latexsym}
\usepackage{amsthm}
\usepackage{amsmath}

\usepackage{amssymb, amscd}

\usepackage{url}

\newtheorem{problem}{Problem}[section]
\newtheorem{theo}[problem]{Theorem}

\newtheorem{defin}[problem]{Definition}
\newtheorem{prop}[problem]{Proposition}

\newtheorem{exam}[problem]{Example}

\begin{document}
\date{May 2, 2006}
 \title{{\Large \bf Groupoids in combinatorics --\\ applications
 of a theory of local symmetries
 }}

\author{Rade  T. \v Zivaljevi\' c\\
Mathematical Institute SANU, Belgrade}

\maketitle

\begin{abstract}
An objective of the theory of {\em combinatorial groupoids} is to
introduce concepts like ``holonomy'', ``parallel transport'',
``bundles'', ``combinatorial curvature'' etc.\ in the context of
simplicial (polyhedral) complexes, posets, graphs, polytopes,
arrangements and other combinatorial objects. In this paper we
give an exposition of some of the currently most active research
themes in this area, offer a unified point of view, and provide a
list of prospective applications in other fields together with a
collection of related open problems.

\end{abstract}

\section{Introduction}

This paper is a sequel to  \cite{Ziv-groupoids} where a program
for developing a theory of combinatorial groupoids was originally
outlined. The main objective of \cite{Ziv-groupoids} was to
demonstrate the relevance of this theory for some well known
problems of contemporary geometric combinatorics, notably the
graph coloring problem (Lov\'asz conjecture and its relatives) and
the problems related to cubifications of manifolds.

In this paper we offer a broader perspective on this subject. More
general concepts are introduced, both old and new applications are
discussed or at least outlined and, what is potentially the most
important aspect of the paper, we try to collect together other
related developments where combinatorial groupoids were implicitly
or explicitly used.

Hoping that this paper may serve as an invitation to the subject,
we included a large number of examples of problems of
combinatorial nature, among them the Penrose impossible ``tribar''
and the S.~Lloyd ``15 game'',  where the ideas and the techniques
of the theory of groupoids may play an important role.

\subsection{An overview}\label{sec:overview}

Recent publications \cite{BoGuHo} \cite{Josw2001}
\cite{Ziv-groupoids} offer a quite convincing evidence that the
language and methods of the theory of groupoids, after being
successfully applied in other major mathematical fields, offer new
insights and perspectives for applications in combinatorics and
discrete and computational geometry.

The groupoids (groups of projectivities) have recently appeared in
geometric combinatorics in the work of M. Joswig \cite{Josw2001},
see also a related paper with Izmestiev \cite{IzmJos2002} and the
references to these papers, where they have been applied to toric
manifolds, branched coverings over $S^3$, colorings of simple
polytopes etc.

\begin{figure}[hbt]
\centering
\includegraphics[scale=0.70]{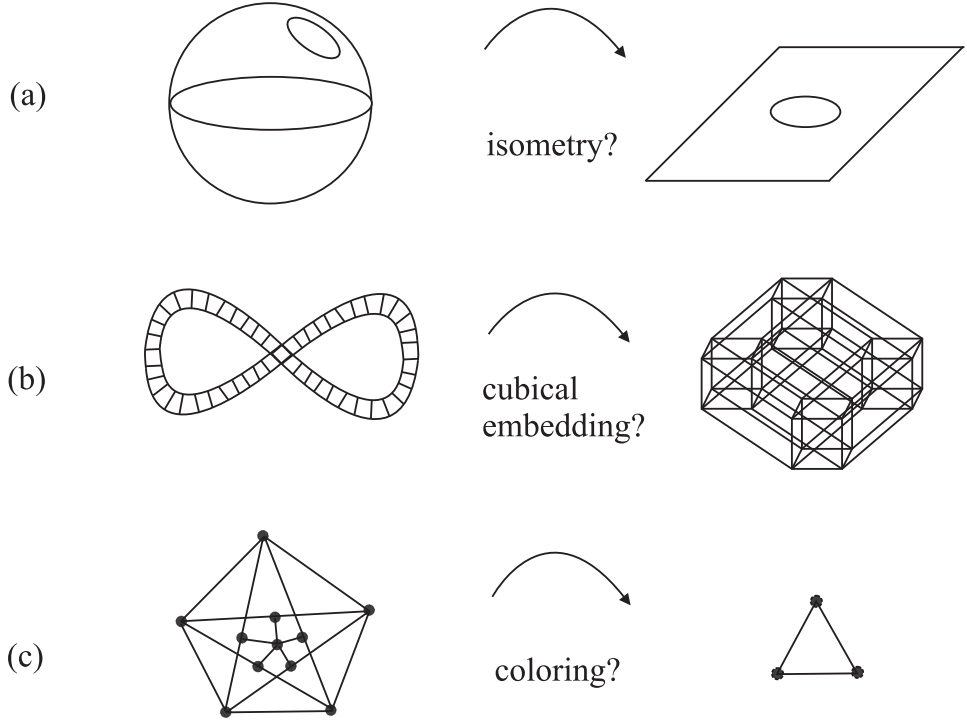}
\caption{Is there a common point of view!?} \label{fig:sphere10}
\end{figure}

Even more explicitly the concepts of ``connection'',
``geodesics'', ``holonomy'' have appeared in \cite{BoGuHo}.
Motivated by the theory of $GKM$-manifolds (named after Goresky,
Kottwitz, and MacPherson \cite{GoKottMcPh}) Bolker, Guillemin, and
Holm develop in this paper an analogy between graph theory and the
theory of manifolds.

The main purpose of \cite{Ziv-groupoids} was to show that these
and related developments should not be seen as isolated examples.
Quite the opposite, they serve as a motivation for further
extensions and generalizations and call for a systematic
applications of a theory of local symmetries in combinatorics.

As a first application it was shown in \cite{Ziv-groupoids} that a
cubical analogue of Joswig's groupoid provides new insight in
cubical complexes non-embedd\-able into cubical lattices (a
question related to a problem of S.P. Novikov which arose in
connection with the $3$-dimensional Ising model) \cite{BuPa02}
\cite{Novikov}. The second, perhaps more far reaching application
developed in this paper, was a generalization, both to more
general test graphs and to simplicial complexes, of a recent
resolution of the Lov\' asz conjecture by Babson and Kozlov
\cite{BabsonKozlov2}.

\begin{figure}[hbt]
\centering
\includegraphics[scale=0.60]{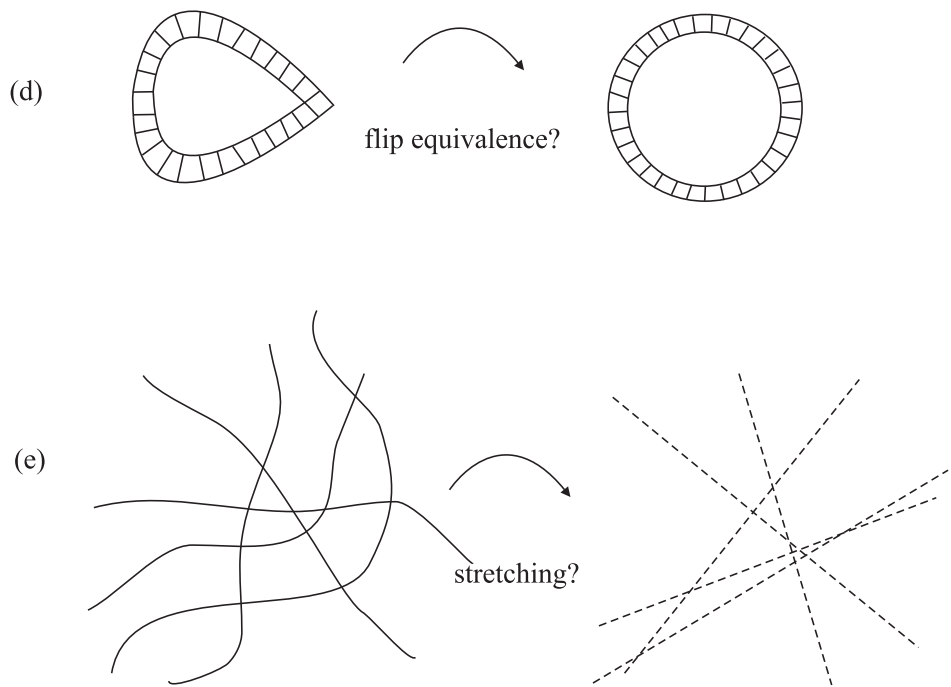}
\caption{Groupoids provide a theory of local symmetries.}
\label{fig:nema-sphere1}
\end{figure}

It appears that combinatorial groupoids are hidden in the
background of many contemporary combinatorial constructions and
applications. It is potentially a very interesting project to
analyze the role of groupoids in the papers like \cite{BBLL}
\cite{BKLW} \cite{BL} \cite{Gaif} etc. One of our central
objectives in this paper is to advocate a systematic use of
groupoids as a valuable tool for geometric and algebraic
combinatorics.

\subsection{The first unifying theme}\label{sec:unifying1}

Our point of departure is an observation that different problems
from different mathematical disciplines, in particular some well
known problems of combinatorial geometric nature, can be all
approached from a similar point of view.

\begin{figure}[hbt]
\centering
\includegraphics[scale=0.50]{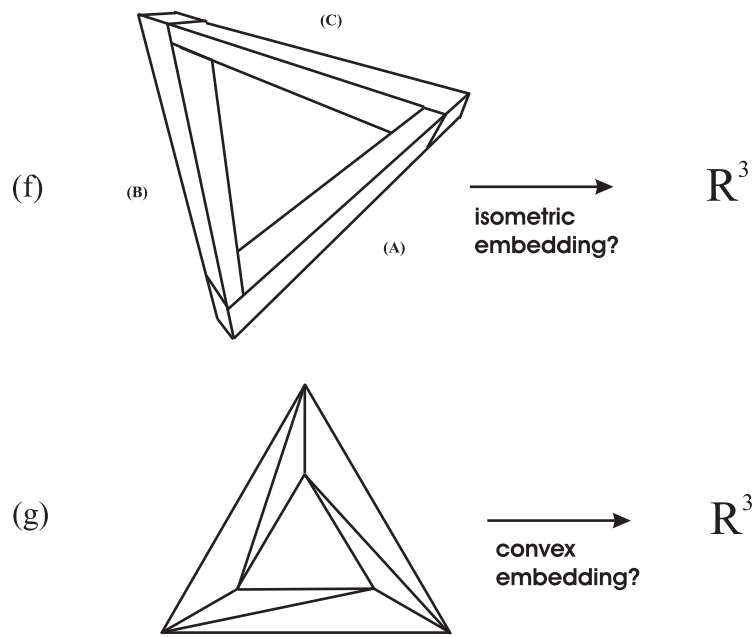}
\caption{Locally possible may be globally impossible!}
\label{fig:tribar6}
\end{figure}

The unifying theme and a single point of view is provided by the
concept of a groupoid. The reader is referred to
Figures~\ref{fig:sphere10}, \ref{fig:nema-sphere1} and
\ref{fig:tribar6} for an informal list of questions which all seem
to involve a concept of a groupoid.

In each of the listed cases there is or ought to be a groupoid
naturally associated to an object of the given category. For each
of these groupoids there is an associated ''parallel transport'',
holonomy groups and other related invariants which serve as
obstructions for the existence of  morphisms indicated in
Figure~\ref{fig:sphere10}.

If $M=(M,g)$ is a Riemannian manifold, the associated groupoid
$\mathcal{G}={\mathcal{G}_{M}}$ has $M$ as the set of objects
while the morphism set $\mathcal{G}(x,y)$ consists of all linear
isomorphisms $a : T_xM\rightarrow T_yM$ arising from a parallel
transport along piece-wise smooth curves from $x$ to $y$. One of
the manifestations of Gauss ``Theorema Egregium'' is that the
associated holonomy group is a metric invariant, consequently an
isometry depicted in Figure~\ref{fig:sphere10}~(a)  is not
possible.

In a similar fashion a cubical ``Theorema Egregium''
\cite[Theorem~3.7]{Ziv-groupoids} provides an obstruction for an
embedding of a cubical complex into a hypercube or a cubical
lattice. As a consequence the cubification (quadrangulation) of
the $\infty$-shaped complex depicted in
Figure~\ref{fig:sphere10}~(b) does not admit a cubical embedding
into a hypercube for the same formal reason (incompatible
holonomy) a spherical cap cannot be isometrically represented in
the plane.

Perhaps it comes as a surprise that the graph coloring problem,
Figure~\ref{fig:sphere10}~(c), can be also approached from a
similar point of view. An analysis of holonomies (``parallel
transport'') of diagrams of $Hom$-complexes of graphs (simplicial
complexes) over the associated Joswig's groupoid
(Example~\ref{ex:Joswig}) eventually leads to a general result
\cite[Theorem~4.21]{Ziv-groupoids} which includes the ``odd'' case
of Babson-Kozlov-Lov\'{a}sz coloring theorem as a special case.

\subsection{The second unifying theme}\label{sec:unifying2}

The problems depicted in Figures~\ref{fig:sphere10},
\ref{fig:nema-sphere1} and \ref{fig:tribar6} can all be seen as
instances of the following general problem-scheme.

\medskip\noindent
{\bf Problem~1:} \cite{Zieg-text} Given some kind of
(combinatorial) structure $\Re$, is it always possible to embed it
into a very ``regular'' or ``complete'' structure of this kind.
Alternatively and perhaps more generally, a ``regular'' or
``complete'' structure may be replaced by some other kind of
environment (space) inhabited by structures similar to $\Re$.

\medskip G.~Ziegler in  \cite{Zieg-text} provides a list of
combinatorial problems which can be placed in this category. They
are all of distinct combinatorial flavor. For example the first
question is whether each matroid of rank $3$ can be embedded into
a finite projective plane, while the second ask if each {\em
Steiner triple system} can be embedded into a finite {\em Kirkman
system} (resolvable Steiner system).

Ziegler makes an interesting remark in this paper to the effect
that ``there should exist cohomology theories that can handle
these embedding problems''. For an interesting evidence that such
``non-classical applications of cohomology theory to embedding
problems'' should exist, he refers the reader to \cite{Pe} and
\cite{CrRy}.

One can speculate that combinatorial groupoids may provide such
``cohomology theory'' in some favorable situations. Indeed, the
essence of the classical Chern-Weyl theory is the construction of
characteristic cohomology classes from the curvature of a manifold
and the curvature is just a manifestation of the holonomy
phenomenon.

For example in the simplest possible situation, an ``obstruction
cocycle'' evaluated on a $2$-dimensional cell, measures the
holonomy around this cell. In other words the information usually
captured by cohomology often comes from the groupoid (connection,
holonomy) naturally associated to the problem.

It is plausible that in majority of "embedding problems" listed in
\cite{Zieg-text}, similarly in each of problems symbolically
depicted in Figures~\ref{fig:sphere10}--\ref{fig:tribar6}, one
should be able to identify combinatorial groupoids which are
naturally associated to these objects (complexes, graphs,
matroids, triple systems, configurations, arrangement etc.).

An embedding induces a morphism of groupoids (often of a very
special type, which implies a monomorphism on the level of
holonomy groups). Already this yields non-embeddability in some
cases (e.g.\ examples of cubical complexes which cannot be
embedded in cubical lattices, Figure~\ref{fig:sphere10} (b)).

However, a cohomology theory we are after ought to be much more
subtle instrument for proving non-embaddibility.

It is fascinating that such a scheme already exists in some sense,
once we identify the groupoids and use objects (matroids, triple
systems, configurations, arrangement) to define natural bundles
over these groupoids.

We pass from groupoids to the associated convolution algebras (in
the same manner one goes from a group to the group algebra or from
a poset to its incidence algebra) and interpret the natural
bundles as moduli over these algebras. After that we are in the
situation which is pregnant with possibilities!

As already emphasized, the classical Chern-Weyl theory is the
construction of characteristic cohomology classes from the
curvature of a manifold, and the relevant information about the
curvature is captured by the underlying groupoid. Today this
construction is incorporated into a map from K-theory of an
algebra (say convolution algebra of a groupoid) to the cyclic
homology of the algebra (Connes, Karoubi etc.).

This is a recipe which in the case of combinatorial groupoids
should be quite concrete and this appears to be a good candidate
for an adequate cohomology theory!

\section{Groupoids}

Groups and symmetries have been treated almost as synonyms in the
history of mathematics. Indeed, we have all been trained that
wherever we encounter symmetries, there ought to be a group of
transformations in the background. Consequently it may come as a
surprise that the concept of a group is sometimes not sufficient
to deal with  this phenomenon in general. Indeed, it may not be
widely known that not groups but {\em grou\-po\-ids} allow us to
handle objects which exhibit what is clearly recognized as
symmetry although they admit no global automorphism whatsoever.
Unlike groups, groupoids are capable of describing reversible
processes which can pass through a number of states. For example
according to A.~Connes \cite{Connes}, Heisenberg discovered
quantum mechanics by considering the groupoid of quantum
transitions rather than the group of symmetry.

\medskip
Groupoids are formally defined as small categories $\mathcal{C} =
(Ob(\mathcal{C}), Mor(\mathcal{C}))$ such that each morphism
$\alpha\in Mor(\mathcal{C}))$ is an isomorphism. This condition
guarantees that each process governed by a groupoid is reversible.
The reader is referred to  \cite{Brown88} \cite{Brown87}
\cite{Brown97} \cite{Higg} \cite{Wein96} for expositions of
different aspects of the theory of groupoids. The {\em vertex}
(isotropy) group $\Pi(\mathcal{C},x):= \mathcal{C}(x,x)$ is often
referred to as the {\em holonomy} group of $\mathcal{C}$ at $x\in
Ob(\mathcal{C})$.
\begin{figure}[hbt]
\centering
\includegraphics[scale=0.50]{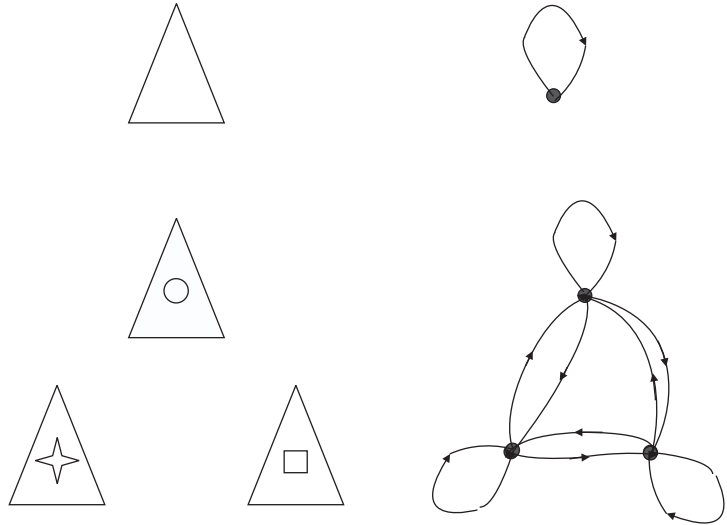}
\caption{From groups to groupoids.} \label{fig:states1}
\end{figure}
A very simple example of a groupoid, which nevertheless exposes
some of its interesting features, is given in
Figure~\ref{fig:states1}. The isosceles triangle has only one
non-trivial symmetry, hence its group of symmetries is
$\mathbb{Z}_2$ . This fact is conveniently recorded in the
directed graph (digraph) depicted on the right of this triangle.
It has only one vertex, corresponding to one object, or one state
of the object, and a directed loop associated to the non-trivial
automorphism.

What if the triangle can pass through a number of different
states, say change its color or if some changing geometric pattern
is present on the triangle? Each ``state'' of the object is
associated a different node in the graph while arrows and directed
loops record all possible transformations between states. The
directed graph obtained this way is precisely the associated
groupoid if we agree that the arrows (transformations) can be
composed, provided they are ``composable'' in the sense that
$\alpha\circ\beta$ exists only if the source object of $\alpha$
coincides with the target object of $\beta$.

\begin{exam}
{\rm An excellent example of how groupoids may appear in
combinatorial practise arises from the analysis of the local
symmetries associated to Penrose tribar,
Figure~\ref{fig:tribar6}~(f). One of the guiding principles of
\cite{Pe} is that this figure is ``locally consistent'' i.e.\ by
covering (removing) one side of the tribar, the remaining two
sides are unambiguously visualized in the surrounding $3$-space.
Moreover, assuming that all three sides are mutually congruent
parallelepipeds, given two of them, say (A) and (B) in
Figure~\ref{fig:tribar6}, there is a natural isometry (local
symmetry) $\alpha_{BA} : (A)\rightarrow (B)$ sending one
parallelepiped to another. In order to make description of
$\alpha_{BA}$ easier to follow, the parallelepipeds are depicted
in Figure~\ref{fig:tribar6} so that one side appears to be thinner
than the other. It is assumed that $\alpha_{BA}$ maps the thinner
side of (A) to the thinner side of (B). The local symmetries
$\alpha_{CB}: (B)\rightarrow (C)$ and $\alpha_{AC}: (C)\rightarrow
(A)$ are similarly defined. It turns out that
 \[
\alpha_{AC}\circ \alpha_{CB}\circ \alpha_{BA} =: \eta_A
 \]
is not an identity map. Rather, it is a rotation of the
parallelepiped (A) through the angle of $90^{\circ}$ around the
longest axes of symmetry.

We conclude that a groupoid associated to the Penrose tribar has
three objects, (A), (B) and (C), and the holonomy group $\Pi\cong
\mathbb{Z}_4$. In other words, the corresponding digraph is very
similar to the digraph depicted in Figure~\ref{fig:states1}.}
\end{exam}

\subsection{Generalities about ``bundles'' and ``parallel transport''}
\label{sec:generalities}

The notion of a groupoid is a common generalization of the
concepts of space and group, i.e.\ the theory of groupoids allows
us to treat spaces, groups and objects associated to them from the
same point of view. In the case of spaces this is achieved by
associating a ``path groupoid'' or the fundamental groupoid $\pi
X$ to a space $X$, \cite[Chapter~6]{Brown88}. Given a group $G$
and a $G$-set $S$, the associated groupoid has $S$ as a set of
objects while morphisms are all ``arrows''
$x\stackrel{g}{\rightarrow} y$ such that $x,y\in S, g\in G$ and
$gx=y$. In particular $G$ itself is a category with only one
object and $G$ as the set of morphisms.

A common generalization for the concept of a bundle $Y$ over $X$
and a $G$-space $Y$ is a $\mathcal{C}$-space or more formally a
diagram over the groupoid $\mathcal{C}$ defined as a functor $F :
\mathcal{C}\rightarrow Top$.

The reader is referred to  \cite{Brown88} \cite{Brown87}
\cite{Brown97} \cite{Wein96} and the references in these sources
for more information about groupoids.

Here we provide only a list (glossary) of some of the basic
concepts associated to groupoids in the form that will allow their
immediate use in subsequent sections.

\bigskip\medskip

\noindent  {\bf GLOSSARY}
\begin{enumerate}
\item[] {\bf Groupoid:}\, A small category $\mathcal{C} = (O,M)$
where $O=Ob(\mathcal{C})$ is the set of objects and
$M=Mor(\mathcal{C})$ the set of morphisms. Informally speaking the
groupoid $\mathcal{C}$ provides a ``road map'' on $O$ which can be
visualized as a digraph as in Figure~\ref{fig:states1}. The vertex
(holonomy) group at $x\in O$ is $\Pi_x:=Hom_{\mathcal{C}}(x,x)$.
If $\mathcal{C}$ is connected all its vertex groups are isomorphic
and often denoted by $\Pi\mathcal{C}$ or $\Pi$.

\item[]{\bf Bundle over $O$:} \, A collection $\mathcal{X} =
\{X_i\}_{i\in O}$ of spaces or sets (fibres) labelled (indexed) by
elements of the set $O$. A bundle often arises from a map $f :
X\rightarrow O$ with $X_i:=f^{-1}(i)$ as the fibre over $i\in O$.

\item[] {\bf Connection on $\mathcal{X}$:}\, A ``connection'' or
``parallel transport'' on the bundle $\mathcal{X}=\{X(i)\}_{i\in
O}$ is a functor (diagram) $\mathcal{P}: \mathcal{C}\rightarrow
Top$ such that $X(i)=\mathcal{P}(i)$ for each $i\in S$. Informally
speaking, the groupoid $\mathcal{C}$ provides a ``road map'' on
$S$, while the functor $\mathcal{P}$ defines the associated
transport from one fibre to another.
\end{enumerate}
It follows from these definitions that a bundle $\mathcal{X}$ is
just a map $O\rightarrow Top$  while a connection extends this map
to a functor $\mathcal{C}\rightarrow Top$.

\subsection{Principal bundles associated to a groupoid}
\label{sec:frame}

There are several ``tautological'' bundles associated to a
groupoid  $\mathcal{C} = (O,M)$. For example one can associate to
$x\in O$ the vertex group $\Pi_x=Hom_{\mathcal{C}}(x,x)$. There is
a natural connection $\mathcal{P} : \mathcal{C}\rightarrow Set$ on
this bundle where for $\alpha\in Hom(x,y), \, \mathcal{P}(\alpha)
: \Pi_x \rightarrow \Pi_y$ is defined by
$\mathcal{P}(\alpha)(\beta):=\alpha\circ\beta\circ\alpha^{-1}$.

\medskip
The concept of a principal or ``frame'' bundle seems to be of
equal importance in applications of (combinatorial) groupoids.
This notion is a natural unification of the concept of a principal
bundle over a topological space and a free $G$-set.

\begin{defin}\label{def:frame}
Suppose that $O$ is a set and assume that $\mathcal{C}=(O,C)$ and
$\mathcal{D}=(O,D)$ are two groupoids on $O$ as the set of
objects. Moreover, assume that $\mathcal{C}$ is a subgroupoid of
$\mathcal{D}$ in the sense that $C\subset D$ and that
$\mathcal{D}$ is connected as a groupoid. Given an object $a\in
O$, define a bundle $Fr = Fr_a : O\rightarrow Set$ by the formula
$Fr_a(y):= \mathcal{D}(a,y)$. This bundle naturally comes with
both a $\mathcal{D}$ and $\mathcal{C}$-connection. The isomorphism
type of this bundle (connection) is independent of the choice of
object $a$ (as a consequence of connectedness of $\mathcal{D}$).
This bundle together with the associated $\mathcal{C}$-connection
is referred to as a $(\mathcal{C},\mathcal{D})$-principal (or
frame) bundle over $O$.
\end{defin}

Usually it is the groupoid $\mathcal{C}$ we are interested in. The
auxiliary groupoid $\mathcal{D}$ often appears as a natural
``ambient'' groupoid for $\mathcal{C}$. For example if
$\mathcal{C}$ is a free $G$-set, then $\mathcal{D}$ is the
groupoid associated to the set $S=Ob(\mathcal{C})$ as a $Q$-set
where $Q\supseteq G$ is the group of all permutations of $S$.

If $\mathcal{C}=\mathcal{G}_M$ is the groupoid associated to a
Riemannian manifold $(M,g)$, described in
Section~\ref{sec:unifying1}, then $\mathcal{D}$ is the groupoid
$Vect_M$ which associates to a pair of points (objects) $(x,y)$ in
$M$, the morphism set $\mathcal{D}(x,y)=Vect(T_xM,T_yM)$ of all
linear isomorphisms from $T_xM$ to $T_yM$. This is the reason why
$Fr= Fr_{\mathcal{C}}$ is also referred to as a frame bundle,
since in this case $Vect(\mathbb{R}^n,T_xM)$ is the set (manifold)
of all $n$-frames in $T_xM$.

This situation arises in all cases where objects of the groupoid
have natural {\em external isomorphisms}, in particular the group
$Hom_{\mathcal{D}}(x,x)$ of external isomorphisms of $x$ may be
larger than $\Pi_x=Hom_{\mathcal{C}}(x,x)$. This is clearly the
case with the groupoid $\mathcal{G}_M$ where the natural group of
symmetries of $T_xM$ is isomorphic to $GL(n,\mathbb{R})$.

All combinatorial groupoids discussed in
Sections~\ref{sec:cg-games} and \ref{sec:comb-group} are of this
kind. In all these examples the natural (external) isomorphisms
are structure preserving bijections associated to these objects.
In the case of the Joswig groupoid $J(K)$, the natural
isomorphisms are bijective simplicial maps of $d$-simplices so the
(external) symmetry group of a $d$-simplex is the group of all
permutations of its vertices, isomorphic to $S_{d+1}$. In the case
of groupoids associated to games the situation is similar. The
external group of symmetries of a position of a game
(Section~\ref{sec:cg-games}) is the group of all permutations of
the pieces, e.g.\ in the ``15 game'' it is the group $S_{15}$. The
external group of symmetries arising in the context of pure
$d$-dimensional, cubical complexes is the group $B_d$ of
symmetries of a $d$-cube etc..

In all these examples there is a tautological ``outer groupoid''
$\mathcal{D}$ and the associated frame bundle $Fr_{\mathcal{C}}$.

\begin{figure}[hbt]
\centering
\includegraphics[scale=0.60]{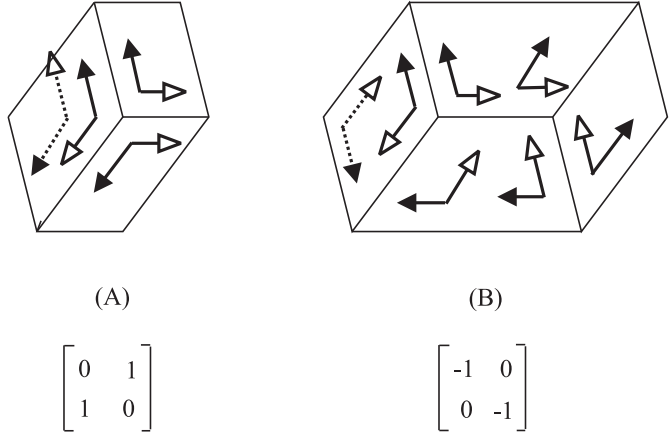}
\caption{Examples of holonomies in a $2$-dimensional cubical
lattice.} \label{fig:dve-putanje1}
\end{figure}

\medskip\noindent
{\bf Symmetry breaking patterns:} Suppose that $\mathcal{C}$ is a
groupoid  where all objects $x\in Ob(\mathcal{C})$ are mutually
{\em externally isomorphic}, i.e.\ isomorphic from the point of
view of their inner (combinatorial or geometric) structure. As a
consequence there is a natural ``external'', connected groupoid
$\mathcal{D}$ associated to $\mathcal{C}$ such that $\mathcal{C}$
is a subgroupoid of $\mathcal{D}$. In other words
$Ob(\mathcal{D})=Ob(\mathcal{C})$ while morphisms in $\mathcal{D}$
are external morphisms. Let $Fr=Fr_a$ be the associated
$(\mathcal{C},\mathcal{D})$-frame bundle. An element of $Fr(x)
=Hom_{\mathcal{D}}(a,x)$ is interpreted as a ``symmetry breaking
pattern'' on $x$. Examples of such patterns are exhibited in
Figures~\ref{fig:trouglici2} and \ref{fig:dve-putanje1} and they
are a useful bookkeeping device for keeping track of the
holonomies, for the concrete combinatorial description of the
associated covering groupoids etc.

\medskip\noindent
{\bf The standard question:} Given $\mathcal{C}$ and the
associated ``outer'' groupoid $\mathcal{D}$, it is interesting to
know whether the associated point groups are different, i.e.\ if
$Hom_{\mathcal{C}}(x,x)\varsubsetneq Hom_{\mathcal{D}}(x,x)$.

\subsection{Combinatorial groupoids; first examples}\label{sec:cg-games}

Combinatorial groupoids are the groupoids that appear in
combinatorics. This is certainly not a very informative statement
so we offer a few examples for illustrative purposes. More formal
definitions are offered in Section~\ref{sec:cg-games}.

\medskip
Suppose that $\mathcal{G}$ is some kind of a ``game'' played on a
``board'' $\mathcal{B}$ with some ``pieces'' $\mathcal{P}$ that
can be moved around this board according to some ``rules''
$\mathcal{R}$. It can be a one-player game, for example a game
with playing cards as pieces (the game of ``Solitaire'' is an
example), a two-player game (chess, checkers etc.) or a
multiplayer game (e.g.\ some multiplayer computer game). We will
ignore here the ``dynamical'' aspect of the game and focus on the
``states'' (positions) of the game $\mathcal{G}$ and how one can,
according to the rules of the game, move from one state to
another.

In order to have a concrete example before our eyes, let us assume
that the board $\mathcal{B}$ is a $(m\times n)$-chessboard and
that pieces cannot be distinguished from one another, like in the
game of checkers. A distribution of pieces on the board
$\mathcal{B}$ is called a position (state) of the game. One can
pass from one position to another by rearranging one or more
pieces, i.e.\ once the game is started pieces are neither removed
from nor returned to the board. An important aspect of this type
of the game is that it is {\em reversible} i.e.\ we can always
return to the original position of pieces by performing the
inverse moves.
\begin{figure}[hbt]
\centering
\includegraphics[scale=0.50]{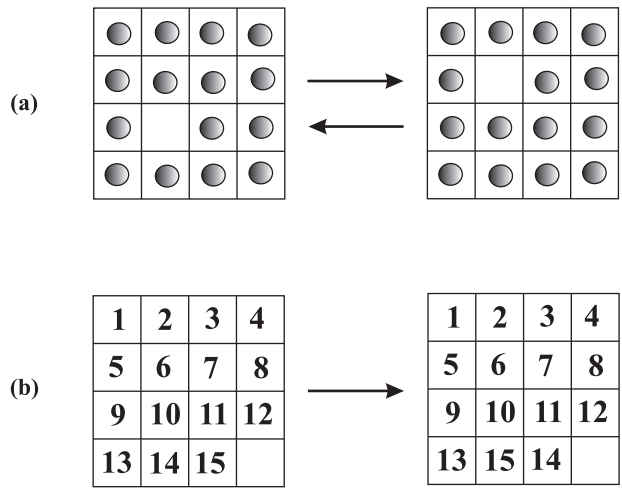}
\caption{The Lloyd ``15 game''.} \label{fig:15game1}
\end{figure}
As customary in game theory, one can associate a directed graph
$D(\mathcal{G})=(N(\mathcal{G}),E(\mathcal{G}))$ to the game
$\mathcal{G}$. The nodes $N(\mathcal{G})$ are all allowed
positions of the game and pairs $(p,q)$ of positions form a
directed edge in $E(\mathcal{G})$ if and only if the rules
$\mathcal{R}$ allow a move from position $p$ to position $q$. It
is clear that this directed graph is actually a groupoid.

It is often convenient to encode all possible positions of the
game in a simplicial complex $K(\mathcal{G})$. The vertices of
$K(\mathcal{G})$ are all elementary cells $(i,j)$ of the board
$\mathcal{B}=[m]\times [n]$ and each position $p\in \mathcal{P}$
of the game contributes a (maximal) simplex $\sigma =\sigma(p)$,
where $(i,j)\in \sigma$ if and only if the cell $(i,j)\in
[m]\times [n]$ is occupied by a piece.

Conversely, given a pure simplicial complex $K$, one can interpret
its maximal simplices as the set of all allowed positions of a
game $\mathcal{G}(K)$ which is played on the set $V(K)$ of
vertices of $K$. Moreover, assume that there is only one rule that
specifies that one can change positions by moving only one piece
at a time. The groupoid arising this way is precisely the Joswig
groupoid $\mathcal{J}(K)$ of $K$ (Section~\ref{sec:Joswig})!

\medskip
One can have even more restrictive rule by asking that only pieces
that satisfy some other constraint can be moved to another
position from a list of allowable positions. A perfect example of
such a game is the famous ``15 game'' of Samuel Lloyd, ``America's
greatest puzzle creator'', see

\noindent
\url{http://www.holotronix.com/samlloyd15a.html}. In this game
$15$ pieces are placed on a $(4\times 4)$-chessboard,
Figure~\ref{fig:15game1}~(a), and a piece can be moved only if it
is an immediate neighbor of the unoccupied cell. Let
$Lloyd_{15}=(O,M)$ be the associated groupoid, i.e.\ the objects
of this groupoid are all $16$ ways to put $15$ identical pieces on
a $(4\times 4)$-chessboard and morphisms are moves allowed by the
``15 game''.

The famous Lloyd's ``15--14'' problem is to start with $15$
labelled pieces, positioned as in Figure~\ref{fig:15game1}~(b) on
the left and,  playing the ``15 game'', end up in the position
depicted in  Figure~\ref{fig:15game1}~(b) on the right. It turns
out that this is not possible. We see this fact as a manifestation
of a phenomenon that
\begin{equation}\label{eqn:holonomy-phenomenon}
\Pi(Lloyd_{15})= A_{15}\varsubsetneqq S_{15}
\end{equation}
i.e.\ that the holonomy group of the groupoid $Lloyd_{15}$ is
different from the a priori given group of symmetries of the
object!

\subsection{Combinatorial groupoids; general picture}
\label{sec:comb-group}

In this section we introduce a sufficiently general class of
combinatorial grou\-po\-ids which seems to capture the essential
features of all examples reviewed in this paper. We warn the
reader that this is certainly not the most general framework
suitable for all possible applications. Rather, as emphasized in
\cite{Ziv-groupoids}, we create ``an ecological niche for
combinatorial groupoids which may be populated by new examples and
variations as the theory develops''.

\begin{defin}\label{def:comb-groupod}
Suppose that $(P,\leq)$ is a (not necessarily finite) poset.
Suppose that $\Sigma$ and $\Delta$ two families of subposets of
$P$. Choose $\sigma_1,\sigma_2\in \Sigma$. If for some $\delta\in
\Delta$ both $\delta\subset \sigma_1$ and $\delta\subset
\sigma_2$, then the posets $\sigma_1,\sigma_2$ are called
$\delta$-adjacent, or just adjacent if $\delta$ is not specified.
Define $\mathcal{C}=(Ob(\mathcal{C}), Mor(\mathcal{C}))$ as a
small category over $Ob(\mathcal{C}) = \Sigma$ as the set of
objects as follows. For two $\delta$-adjacent objects $\sigma_1$
and $\sigma_2$, an {\em elementary morphism} $\alpha\in
\mathcal{C}(\sigma_1,\sigma_2)$ is an isomorphism $\alpha :
\sigma_1\rightarrow \sigma_2$ of posets which leaves $\delta$
point-wise fixed. A morphism $\mathfrak{p} \in
\mathcal{C}(\sigma_0,\sigma_m)$ from $\sigma_0$ to $\sigma_m$ is
an isomorphism of posets $\sigma_0$ and $\sigma_m$ which can be
expressed as a composition of elementary morphisms.
\end{defin}

Given two adjacent objects $\sigma_1$ and $\sigma_2$, an
elementary morphism $\alpha\in \mathcal{C}(\sigma_1,\sigma_2)$ may
not exist at all, or if it exists it may not be unique. In case it
exists and is unique it will be frequently denoted by
$\overrightarrow{\sigma_1\sigma_2}$ and sometimes referred to as a
``flip'' from $\sigma_1$ to $\sigma_2$. In this case a morphism
$\mathfrak{p}\in \mathcal{C}(\sigma_0,\sigma_m)$ is by definition
a composition of flips
\begin{equation}\label{eqn:flips}
\mathfrak{p} = \overrightarrow{\sigma_0\sigma_1}\ast
\overrightarrow{\sigma_1\sigma_2}\ast\ldots\ast
\overrightarrow{\sigma_{n-1}\sigma_n}.
\end{equation}

\medskip\noindent
{\bf Caveat:} Here we adopt a useful convention that $(x)(f\ast
g)= (g\circ f)(x)$ for each two composable maps $f$ and $g$. The
notation $f\ast g$ is often given priority over the usual $g\circ
f$ if we want to emphasize that the functions act on the points
from the right, that is if the arrows in the associated formulas
point from left to the right.

\medskip

Suppose that $P$ is a ranked poset of depth $n$ with the
associated rank function $r : P\rightarrow [n]$. Let $\mathcal{E}
= \mathcal{E}_P$ be the $\mathcal{C}$-groupoid described in
Definition~\ref{def:comb-groupod} associated to the families
$\Sigma :=\{P_{\leq x} \mid r(x)=n\}$ and $\Delta :=\{P_{\leq y}
\mid r(y)=n-1\}$. It is clear that other ``rank selected''
groupoids can be similarly defined.

\medskip
The definitions of groupoids $\mathcal{C}$ and $\mathcal{E}$ are
easily extended from posets to simplicial, polyhedral, or other
classes of cell complexes. If $K$ is a complex and $P := P_K$ the
associated face poset, then $\mathcal{C}_K$ and
$\mathcal{E}_K=\mathcal{E}(K)$ are groupoids associated to the
poset $P_K$. We will usually drop the subscript whenever it is
clear from the context what is the ambient poset $P$ or complex
$K$.

\begin{exam}\label{ex:Joswig}{\rm Suppose that $K$ is a pure,
$d$-dimensional simplicial complex. Let $\mathcal{E}(K)$ be the
associated $\mathcal{E}$-groupoid corresponding to ranks $d$ and
$d-1$. Then the {\em groups of projectivities} $\Pi(K,\sigma)$,
introduced by Joswig in \cite{Josw2001}, are nothing but the
holonomy groups of the groupoid $\mathcal{E}(K)$. For this reason
the groupoid $\mathcal{E}(K)$ is in the sequel often referred to
as Joswig's groupoid and denoted by $\mathcal{J}(K)$.
$\mathcal{J}(K)$ is {\em connected} as a groupoid if and only if
$K$ is ``strongly connected'' in the sense of \cite{Josw2001}.
 }
\end{exam}

A simplicial map of simplicial complexes is non-degenerate if it
is 1--1 on simplices. The following definition extends this
concept to the case of posets.

\begin{defin}\label{def:nondeg}
A monotone map of posets $f : P\rightarrow Q$ is {non-degenerate}
if the restriction of $f$ on $P_{\leq x}$ induces an isomorphism
of posets $P_{\leq x}$ and $Q_{\leq f(x)}$ for each element $x\in
P$. Similarly, a map of simplicial, cubical or more general cell
complexes is non-degenerate if the associated map of posets is
non-degenerate. In this case we say that $P$ is {\em mappable} to
$Q$ while a non-degenerate map $f : P\rightarrow Q$ is often
referred to as a {\em combinatorial immersion} from $P$ to $Q$.
\end{defin}

\begin{exam}\label{ex:graph-homomorphism}
{\rm A graph homomorphism \cite{Kozlov-review} $f : G_1
\rightarrow G_2$ can be defined as a non-degenerate map of
associated $1$-dimensional cell complexes. A $n$-coloring of a
graph $G$ is a non-degenerate map (graph homomorphism) $f :
G\rightarrow K_n$ where $K_n$ is a complete graph on $n$
vertices.}
\end{exam}

\begin{prop}\label{prop:functor}
Suppose that $P$ and $Q$ are ranked posets of depth $n$ and let $f
: P\rightarrow Q$ be a non-degenerate map. Then there is an
induced map (functor) $F : \mathcal{E}_P \rightarrow
\mathcal{E}_Q$ of the associated $\mathcal{E}$-groupoids.
Moreover, $F$ induces an inclusion map $\Pi(\mathcal{E}_P,
\mathfrak{p})\hookrightarrow \Pi(\mathcal{E}_Q, F(\mathfrak{p}))$
of the associated holonomy groups.
\end{prop}

\section{Applications}\label{sec:applications}

\subsection{Joswig groupoid $\mathcal{J}(K)$}\label{sec:Joswig}

M.~Joswig defined parallel transport and the associated holonomy
groups in the context of each pure, $d$-dimensional simplicial
complex $K$ in \cite{Josw2001}. He did not formally use the
language of the theory of groupoids, but the combinatorial
groupoids are implicit in this and subsequent paper
\cite{Jos-Usp}, and in the joint paper with I.~Izmestiev
\cite{IzmJos2002}. In particular, our definition of combinatorial
groupoids and associated concepts is strongly influenced by
Joswig's point of view and reflects a desire to incorporate other
examples of apparently similar nature into the same framework.

Formally speaking, the Joswig groupoid is the
$\mathcal{E}$-groupoid (Section~\ref{sec:comb-group}) associated
to a pure $d$-dimensional simplicial complex $K$
(Example~\ref{ex:Joswig}).

\begin{figure}[hbt]
\centering
\includegraphics[scale=0.50]{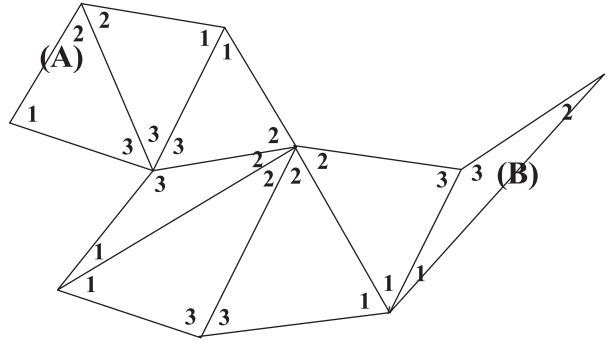}
\caption{Parallel transport from \bf{(A)} to \bf{(B)}.}
\label{fig:trouglici2}
\end{figure}

More explicitly the objects of $\mathcal{K}$ are $d$-dimensional
simplices of $K$ while the morphisms are compositions of ``flips''
(Figure~\ref{fig:trouglici2}), as in the equation
(\ref{eqn:flips}), Section~\ref{sec:comb-group}.

A ``symmetry braking pattern'', in the sense of
Section~\ref{sec:frame}, which can be used for keeping the track
of the action of a Joswig groupoid is simply a labelling of
vertices of the corresponding simplex, see e.g.\ simplex~(A) in
Figure~\ref{fig:trouglici2}. The holonomy group is a priory a
subgroup of the group of all permutations of vertices of a
simplex, i.e.\ a subgroup of the group $S_{d+1}$.

\medskip
Let us briefly review the original problem that motivated
M.~Joswig to introduce his ``groups of projectivities'', i.e.\ the
holonomy groups of the groupoid $\mathcal{J}(K)$.

M.~Davis and T.~Januszkiewicz associated a smooth
$(n+d)$-dimensional manifold $\mathcal{Z}_P$ to any
$d$-dimensional, simple convex polytope $P$ with $n$ facets. These
are examples of quasi-toric manifolds \cite{BuPa02}, relatives of
toric varieties, which come equipped with a $T^n$-action. Genuine
toric varieties are of dimension $2d\leq n+d$. Perhaps motivated
in part by this, V.~Buchstaber suggested a program of studying
when one can find a subgroup $T$ of $T^n$ which acts freely on
$\mathcal{Z}_P$; in that case $\mathcal{Z}_P/T$ would be another
quasi-toric manifold of dimension $n+d-{\rm dim}(T)$. Let $s(P)$
is the maximum dimension of a subgroup $T\subseteq T^n$, acting
freely on $\mathcal{Z}_P$. I.~Izmestiev in \cite{Izm2001} defined
the chromatic number $\gamma(P)$ of $P$ as the minimal number of
colors required to color the facets of $P$ such that any two
facets sharing a vertex have different colors. The relation
between $s(P)$ and $\gamma(P)$ is given by inequalities
\[
n-\gamma(P)\leq s(P)\leq n-d
\]
where the right hand inequality is elementary while the left hand
relation is due to Izmestiev \cite{Izm2001}.

One of the consequences of the Joswig's  analysis of holonomy
groups of the groupoid $\mathcal{J}(K)$ where $K$ is the dual of
$P$, is the result \cite[Theorem~16]{Josw2001} which implies that
$\gamma(P)=d$ if the corresponding holonomy group is trivial.

\medskip
The reader is referred to \cite{Izm2001} \cite{Izm-Usp2001}
\cite{IzmJos2002}   \cite{Josw2001} \cite{Jos-Usp} for these and
other applications of group of projectivities of simplicial
complexes.

\subsection{How is graph like a manifold?}\label{sec:graph=manifold}

One of the central ideas of \cite{BoGuHo} is to approach classical
combinatorial problems by exploring analogies between graphs and
manifolds. A central theme, illuminating this connection, arises
from the theory of group actions, notably from the analysis of
complex $(\mathbb{{C}^\ast})^n$-manifolds. The class of so called
$GKM$-manifolds, named after Goresky, Kottwitz, and MacPherson,
has a particularly interesting structure theory. More precisely,
the $0$-dimen\-sio\-nal orbits, as nodes, and $1$-dimensional
orbits, as edges, define an associated $GKM$-graph $\Gamma=
\Gamma(M)$ which captures a substantial part of the structure of
the original $GKM$-manifold $M$. The graph $\Gamma$ arising this
way turns out to be $d$-regular, where $d$ is the dimension of the
underlying complex $\mathbb{C}^n$-manifold $M$. An extra piece of
structure is an assignment of integer vectors (axial functions) to
edges of this graph, which taken together define an ``embedding''
of the graph in $\mathbb{R}^n$.

\medskip
Bolker, Guillemin, and Holm, building on the previous work of
Goresky, Kottwitz, MacPherson, Rosu, Knutson, Lian, Liu, Yau,
Zara, and others, develop in this paper a dictionary associating
manifold concepts to graph concepts.

A particularly interesting aspect of this work is appearance of
{\em connections, holonomy groups, geodesics}, etc.\ in the
context of arbitrary (regular) graphs $\Gamma = (V,N)$\footnote{I
am grateful to M.~Joswig for drawing my attention to this fact!}.
Here is one of the main definitions.

\begin{defin}{\rm \cite{BoGuHo}} A {\em connection} on a graph
$\Gamma = (V,E)$ is a collection of bijective functions
$\nabla_{(x,y)}: Star(x)\rightarrow Star(y)$, indexed by all
(oriented!) edges $(x,y)$ in $\Gamma$, where $Star(z):=\{(z,w)\in
E\mid w\in V\}$ is the set all oriented edges in $\Gamma$ incident
to $z$.  These functions satisfy the following conditions:

\begin{enumerate}
 \item[{\rm (1)}] $\nabla_{(x,y)}(x,y)=(y,x)$,
 \item[{\rm (2)}] $\nabla_{(y,x)}=\nabla_{(x,y)}^{-1}$.
\end{enumerate}
\end{defin}
Bolker, Guillemin, and Holm use the connection $\nabla$ to define
geodesics in the graph $\Gamma$, to introduce its totally geodesic
subgraphs, holonomy groups as subgroups of the groups of all
permutations of $Star(x)$ etc., see \cite{BoGuHo} for the detailed
development and applications of these concepts.

The reader is invited to identify the associated combinatorial
groupoid $\mathcal{G}_{\Gamma}$ and to relate it to the groupoids
described in Sections~\ref{sec:cg-games} and \ref{sec:comb-group}.
The associated bundle where this connection (parallel transport)
is defined) is clearly the collection $\{Star(x)\}_{x\in V}$.

Answering the question from the title of their paper (and our
Section~\ref{sec:graph=manifold}), Bolker, Guillemin, and Holm in
\cite[Section~2]{BoGuHo} state that:``... the star of a vertex (of
a graph) is a combinatorial analogue of the tangent space to a
manifold at a point ...''.

We observe that this is in complete agreement with the point of
view of our introductory sections. Indeed, a ``tangent space'' is
in all exhibited examples either an object in $Ob(\mathcal{G})$ or
alternatively the fibre of a tautological bundle over the
associated groupoid $\mathcal{G}$.

\subsection{Holonomy vs.\ NaCl-invariant of a cubical complex}
\label{sec:NaCl}

In this section we apply the ideas outlined in earlier sections to
the case of cubical complexes.

Recall that a cell complex $K$ is cubical if  it is a regular
$CW$-complex such that the associated face poset $P_K$ is cubical
in the sense of the following definition.

\begin{defin}\label{def:cubical}
$P$ is a cubical poset if:
\begin{enumerate}
 \item[{\rm (a)}] for each $x\in P$, the subposet $P_{\leq x}$ is
isomorphic to the face poset of some cube $I^q;$
 \item[{\rm (b)}] $P$ is a semilattice in the sense that if a pair
 $x,y\in P$ is bounded from above then it has the least upper
 bound.
\end{enumerate}
\end{defin}

If a space $X$ comes equipped with a standard cubification, clear
from the context, this cubical complex is denoted by $\{X\}$, the
associated $k$-skeleton is denoted by $\{X\}_{(k)}$ etc. For
example $\{I^d\}_{(k)}$ is the $k$-skeleton of the standard
cubification of the $d$-cube, similarly $\{\mathbb{R}^d\}_{(k)}$
is the $k$-skeleton of the standard cubification of $\mathbb{R}^d$
associated to the lattice $\mathbb{Z}^d$.

\begin{figure}[hbt]
\centering
\includegraphics[scale=0.50]{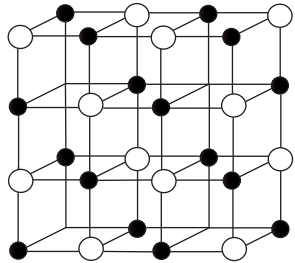}
\caption{Sodium chloride NaCl as a cubical complex.}
\label{fig:NaCl1}
\end{figure}

The group $B_k$ of all symmetries of a $k$-cube is isomorphic to
the group of all signed, permutation $(k\times k)$-matrices. Its
subgroup of all matrices with even number of $(-1)$-entries is
denoted by $B_k^{even}$. The vertex-edge graph of a $k$-cube is
well known to be bipartite, i.e.\ colorable with two colors
(Figure~\ref{fig:NaCl1}) and $B_k^{even}$ can be described as the
set of all elements in $B_k$ that preserve this coloring.

Given a (pure) $d$-dimensional cubical complex $X$, the associated
groupoid is denoted by $\mathcal{C}(X)$. As a cubical counterpart
of Joswig's groupoid, it was introduced in \cite{Ziv2005}
\cite{Ziv-groupoids} and applied to problems related to embeddings
of cubical complexes into cubical lattices (problem of S.~Novikov,
Figure~\ref{fig:sphere10} (b)) and questions of ``bubble
modifications'' of cubical complexes (problem of N.~Habbeger,
Figure~\ref{fig:nema-sphere1} (d)).

Both applications were based on a holonomy type,
$\mathbb{Z}_2$-invariant $I(K)$ of a cubical complex $K$
introduced in \cite{Ziv-groupoids}.

\begin{defin}
Suppose that $K$ is a $k$-dimensional cubical complex and let
$\Pi(K,\sigma)$ be its combinatorial holonomy group based at
$\sigma\in K$. By definition let $I(K)= 0$ if
$\Pi(K,\sigma)\subset B_k^{even}$ for all $\sigma$, and $I(K)= 1$
in the opposite case.
\end{defin}

The reader is referred to \cite{Ziv-groupoids} for more detailed
exposition of results related to this invariant. Following a
suggestion of G.~Ziegler, we give a useful criterion which in many
cases of interest enables us to prove that $I(K)=0$.

\begin{prop}\label{prop:NaCl}(NaCl-criterion) $I(K)=0$ if the cubical
complex $K$ can be colored with two colors such that adjacent
vertices are always of different colors, equivalently if the
vertex-edge graph of $K$ is bipartite.
\end{prop}

\medskip\noindent
{\bf Proof:} (outline) If such a coloring exists then it is
preserved by the parallel transport in the groupoid
$\mathcal{C}(K)$. \hfill$\square$

Let us define the ``NaCl-invariant'' of a cubical complex $K$ by
the requirement that ${\rm NaCl}(K)$ is $0$ (respectively $1$) if
its vertex-edge graph can (cannot) be colored by $2$ colors
(Figure~\ref{fig:NaCl1}). Proposition~\ref{prop:NaCl} says that
$I(K)\leq {\rm NaCl}(K)$. It is not difficult to find examples of
complexes which prove that in general $I(K)\neq {\rm NaCl}(K)$.

\begin{exam}{\rm
Indeed, let $K$ be a cubical complex such that ${\rm NaCl}(K)=0$.
Let us identify two vertices $u$ and $v$, non-adjacent in $K$,
which are nevertheless assigned different colors. Let $K'$ be the
cubical complex $K':= K/u\thickapprox v$. This identification does
not effect the holonomy group of $K$ i.e.\ $\Pi_K\cong
\Pi_{K'}\subset B_k^{even}$, hence $I(K)=I(K')=0$. On the other
hand ${\rm NaCl}(K')=1$.}
\end{exam}

Let us clarify the relationship between invariants $I(K)$ and
${\rm NaCl}(K)$, at least in the important case of cubifications
of manifolds. As it was kindly pointed by the referee, the
following proposition, modelled on Proposition~6 from
\cite{Josw2001}, provides a fairly complete and natural answer to
this question.

\begin{prop}
Suppose that $K$ is a cubical complex which is globally and
locally (strongly) connected. This means that both the groupoid
$\mathcal{C}(K)$ and each of its subgroupoids $\mathcal{C}({\rm
Star}(v))$ are connected where $v$ is a vertex in $K$. Then
\[
I(K) = {\rm NaCl}(K).
\]
\end{prop}

\medskip\noindent
{\bf Proof:} The proof is similar to the proof of Proposition~6
from \cite{Josw2001}. Since $I(K)\leq {\rm NaCl}(K)$, it is
sufficient to show that if $I(K)=0$ then the vertices of $K$ can
be colored with two colors such that no two vertices with same
color are adjacent, i.e.\ that the vertex-edge graph of $K$ is
bipartite.

The required coloring of vertices of $K$ is defined as follows.
Select a top dimensional cell $c_0\in K$ and color its vertices
with two colors. Each other cell $c\in K$ is (strongly) connected
with $c_0$, i.e.\ connected in the sense of the groupoid
$\mathcal{C}(K)$. Given a ``path'' (morphism) $\theta$ in
$\mathcal{C}(K)$, connecting $c_0$ and $c$, the coloring of $c_0$
can be extended along this path to a coloring of $c$. By
assumption $\Pi(\mathcal{C}(K))=0$, hence this coloring of $c$
does not depend of the path $\theta$. On the other hand a vertex
$v\in c$ might receive a different color from another cell $c'$
such that $v\in c'$. By assumption ${\rm Star}(v)$ is also
(strongly) connected which guaranties that there exists a path
from $c$ to $c'$ inside ${\rm Star}(v)$. This guarantees that this
is not possible which concludes the proof of the proposition.
\hfill$\square$

\medskip
The following result, illustrates the elegance of the ${\rm
NaCl}$-approach. It shows that Theorem~3.2 from
\cite{Ziv-groupoids} is really a relative of results from
\cite{Josw2001} about colorings of simple polytopes.

\begin{theo}\label{thm:cubical-lattice}
Suppose that $K$ is a $k$-dimensional cubical complex which is
embeddable/mappable to $\{Z\}_{(k)}$, the $k$-dimensional skeleton
of the standard cubical decomposition of a generic zonotope
$Z=[-v_1,v_1]+\ldots +[-v_n,v_n]$. Then $I(K)=0$.
\end{theo}

\medskip\noindent
{\bf Proof:} Each vertex $w = \epsilon_1v_1+\ldots +\epsilon_nv_n$
of $Z$ is colored in ``black'' (respectively'``white'') if there
is an odd (even) number of occurrences of $-1$ in the sequence
$\epsilon_1,\ldots,\epsilon_n$. It is not difficult to show that
in this coloring no two adjacent vertices are colored by the same
color. This implies that ${\rm NaCl}(K)=0$ for each
$k$-dimensional subcomplex of $Z$.\hfill $\square$

\subsection{Generalized Lov\'{a}sz conjecture}

One of the central problems of topological graph theory is to
explore how the topological complexity of a graph complex $X(G)$
is reflected in the combinatorial complexity of the graph $G$
itself. The results one is often interested in come in the form of
implications
\[
\alpha(X(G))\geq p \Rightarrow \xi(G)\geq q,
\]
where $\alpha(X(G))$ is a topological invariant of the complex
$X(G)$, while $\xi(G)$ is a combinatorial invariant of the graph
$G$. The earliest statement of this type is the celebrated result
of L.~Lov\' asz which is today often formulated in the form of an
implication
\[
Hom(K_2,G) \mbox{ {\rm is $k$-connected} } \Rightarrow \chi(G)\geq
k+3,
 \]
where $Hom(K_2,G)$ is one of many (essentially equivalent)
$\mathbb{Z}_2$-complexes associated to $G$, see
\cite{Kozlov-review} and \cite{Matousek} as overviews and guides
to the literature. $Hom(K_2,G)$ is a special case of a general
graph complex $Hom(H,G)$ (also introduced by L.~Lov\' asz), a cell
complex which functorially depends on the input graphs $H$ and
$G$.

An outstanding conjecture in this area, referred to as the ``Lov\'
asz conjecture'', was that one obtains a better bound if the graph
$K_2$ is replaced by an odd cycle $C_{2r+1}$. More precisely Lov\'
asz conjectured that
\begin{equation}\label{eq:LC}
Hom(C_{2r+1},G) \mbox{ {\rm is $k$-connected} } \Rightarrow
\chi(G)\geq k+4.
 \end{equation}
This conjecture was confirmed in \cite{BabsonKozlov2}, see also
\cite{Kozlov-review} for a more detailed account.

\medskip
The main observation of \cite{Ziv2005} was that the general theory
of groupoids, in particular the Joswig groupoid $\mathcal{J}(K)$,
provide a deep insight into the L\'{o}vasz conjecture and its
ramifications. As a consequence one obtains the implication
\[
Hom(\Gamma,K) \mbox{ {\rm is $k$-connected} } \Rightarrow
\chi(K)\geq k+d+3
\]
which under  suitable assumption on the ``test complex'' $\Gamma$
and the assumption that integer $k$ is odd, extends the result of
Babson and Kozlov to the case of pure $d$-dimensional simplicial
complexes. Moreover, this approach yields a short and conceptual
proof of the Lov\' asz conjecture for $k$ odd. The reader is
referred to \cite{Ziv-groupoids} for an exposition of these and
related results.

Subsequently the approach based on groupoids was extended and
incorporated into equivariant index theory by C.~Schultz
\cite{Schultz} \cite{Schultz06} who developed new powerful methods
leading to deep understanding of $Hom$-complexes and further
analogues of (\ref{eq:LC}).

There are two new, short and elegant, proofs of the
Babson-Kozlov-Lov\'{a}sz theorem. The proof in \cite{Schultz06} is
based on the evaluation of the cohomological $\mathbb{Z}_2$-index
while the more recent proof \cite{Kozlov-short} relies on a
combinatorial evaluation of the height of the associated
Stiefel-Whitney characteristic class.

\subsubsection{The main observation}

In this section we briefly describe the nature of the
``mathematical revelation'' that pointed to the connection between
\cite{BabsonKozlov2} and \cite{Josw2001}, led to \cite{Ziv2005}
and \cite{Ziv-groupoids}, and served as the author's main initial
motivation for starting the program of studying combinatorial
groupoids.

\medskip
It is well known that a graph $G=(V_G,E_G)$ admits a coloring with
not more than $m$ colors if and only there exists a graph
homomorphism $c: G\rightarrow K_m$ from $G$ to the complete graph
with $m$ vertices (Example~\ref{ex:graph-homomorphism}).

Given graphs $G=(V_G,E_G)$ and $H=(V_H,E_H)$, the associated graph
complex $Hom(G,H)$ is the cell complex where each cell is indexed
by a multivalued function $\eta : V_G\rightarrow
2^{V_H}\setminus\{\emptyset\}$ such that if $(i,j)\in E_G$ then
for each $\alpha\in\eta(i)$ and each $\beta\in\eta(j),\,
(\alpha,\beta)\in E_H$, \cite{Kozlov-review} \cite{BabsonKozlov2}
\cite{Schultz06}. For example it is a well known fact that the
$Hom$-complex $Hom(K_2,K_m)$ between complete graphs $K_2$ and
$K_m$ is homeomorphic to a $(m-2)$-dimensional sphere.

\medskip
A graph $G$ (without loops and multiple edges) can be interpreted
as a $1$-dimensional simplicial complex. Let $\mathcal{J}(G)$ be
the corresponding Joswig's groupoid, Example~\ref{ex:Joswig}. Each
edge $e\in E_G$ itself can be interpreted as a subgraph of $G$
isomorphic to $K_2$. The map $\Gamma_G : E_G \rightarrow Top$
defined by
\[
e\mapsto Hom(e,K_m)
\]
is a spherical bundle over the set $E_G$ of edges of $G$ in the
sense of Section~\ref{sec:generalities}. There is a ``forgetful''
continuous map $\phi_e: Hom(G,K_m)\rightarrow Hom(e,K_m)$ for each
edge $e\in E_G$.

\medskip\noindent
{\bf The key observation}: The ``parallel transport'' with respect
to the Joswig's groupoid $\mathcal{J}(G)$ preserves the homotopy
type of the map $\phi_e$. If $G\cong C_{2r+1}$ is an odd cycle,
then the holonomy group $\Pi(\mathcal{J}(G))\cong \mathbb{Z}_2$ is
nontrivial and as a consequence there is a homotopy equivalence
\begin{equation}\label{eq:homotopy}
\phi_e\simeq \phi_e\circ \alpha_e
\end{equation}
where $\alpha_e\in \Pi_e$ is the nontrivial element of the
corresponding holonomy group.

The homotopy (\ref{eq:homotopy}) has cohomological consequences
which eventually, in light of the naturality of
$Hom$-construction, can be used to show that a coloring
$c:G\rightarrow K_m$ is not possible. The details of this
construction and its ramifications can be found in
\cite{Ziv-groupoids}, see also \cite{Ziv2005} for a preliminary
version.

\subsection{Afterword}

There are  other examples of applications of discrete connections,
discrete holonomies (combinatorial groupoids) etc.\ that have not
been covered by this review. A notable example is the paper of
Novikov \cite{Novikov04}, see also Novikov and Dynnikov
\cite{NovDyn}, and the references in these papers.

Novikov and his followers have studied discrete connections on
triangulated manifolds as a part of a general programme of
developing discretized differential geometry, finding discrete
analogs of important differential operators, describing discrete
analogs of complete integrable systems etc.

These developments are naturally linked with the ``Discrete
differential geometry'' in the sense of Bobenko and Suris
\cite{BobSur}, a broad new area where differential geometry of
smooth curves, surfaces and other manifolds interacts with
discrete geometry, using tools and ideas from all parts of
mathematics, and having applications ranging from integrable
dynamical systems to computer graphics.

\bigskip
\noindent {\bf\large Acknowledgement:} It is a pleasure to
acknowledge encouragement, useful comments and suggestions by
R.~Brown, M.~Joswig, J.~Sullivan, G.~Ziegler, the referee, and
numerous participants of the conferences ``Algebraic and Geometric
Combinatorics'', Anogia (Crete), August 20--26, 2005, and
``Discrete Differential Geometry'', Mathematisches
Forschungsinstitut Oberwolfach, March 5--11, 2006.

This is also a pleasant opportunity to acknowledge the support by
the projects no.\ 144014 and 144026 of the Serbian Ministry of
Science and the project ``Geometry, Topology and Combinatorics of
Manifolds and Dynamical Systems'' (SISSA, Trieste), of the Italian
Ministry of Universities and Scientific Research.


\vfill\newpage
 \small \baselineskip3pt
\begin{thebibliography}{abcdefg}


\bibitem[BBLL]{BBLL} E.~Babson, H.~Barcelo, M. De Longueville,
and R.~Laubenbacher, \textit{Homotopy theory for graphs},
arXiv:math.CO/0403146 v1 9 Mar 2004, Journal of Algebraic
Combinatorics, in press.

\bibitem[BC]{BC} E.K.~Babson and C.~Chan, \textit{Counting faces for cubical
spheres modulo two}, Discrete Math.\ 212 (2000), 169--183.

\bibitem[BK03]{BabsonKozlov1} E.~Babson, D.N.~Kozlov,
\textit{Complexes of graph homomorphisms}, \newline
arXiv:math.CO/0310056 v1 5 Oct 2003, to appear in Israel J.\ Math.

\bibitem[BK04]{BabsonKozlov2} E.~Babson, D.N.~Kozlov,
\textit{Proof of the Lov\' asz conjecture}, Annals of Mathematics,
Accepted papers (2005), arXiv:math.CO/0402395 v2, 2004.


\bibitem[BL]{BL} H.~Barcelo and R.~Laubenbacher, \textit{Perspectives on
$A$-homotopy theory and its applications}, preprint.

\bibitem[BKLW]{BKLW} H.~Barcelo, X.~Kramer, R.~Laubenbacher, and C.~Weaver,
\textit{Faundations of a connectivity theory for simplicial
complexes}, preprint.

\bibitem[BEE]{BEE} M.W.~Bern, D.~Eppstein, and J.G.\ Erikson,
\textit{Flipping cubical meshes}, Engineering with Computers 18
(2002), 173--187.

\bibitem[Bj95]{Bjorner} A. Bj\"{o}rner, \textit{Topological methods}.
In R. Graham, M.\ Gr\"{o}tschel, and L. Lov\'{a}sz, editors,
\textit{Handbook of Combinatorics}, North-Holland, Amsterdam,
1995.


\bibitem[BVSWZ]{Bjo-et-al} A.~Bj\" orner, M.~Las Vergnas, B.~Sturmfels,
N.~White, G.~Ziegler, {\it Oriented Matroids}, Encyclopedia of
Mathematics and its Applications \textbf{46}, Cambridge University
Press 1993.


\bibitem[BS]{BobSur} A.I.~Bobenko, Yu.B.~Suris, \textit{Discrete
Differential Geometry; Consistency as Integrability},
arXiv:math.DG/0504358 v1 Apr 2005.


\bibitem[BGH]{BoGuHo} E.D.~Bolker, V.W.~Guillemin, T.S.~Holm,
\textit{How is a graph like a manifold?}, arXiv:math.CO/0206103 v1
Jun 2002.


\bibitem[Br87]{Brown87} R.~Brown, \textit{From groups to groupoids: a
brief surway}, Bull.\ London Math.\ Soc.\ 19(1987) 113--134.


\bibitem[Br97]{Brown97} R.~Brown, \textit{Groupoids and crossed objects in
algebraic topology}, Homology, Homotopy and Applications, Vol.\ 1,
1999, No.\ 1, pp 1--78.

\bibitem[Br]{Brown88} R.~Brown, \textit{Topology and Groupoids}, Booksurge LLC
2006; retitled,  revised, updated and extended edition of
\textit{Topology: A Geometric Account of General Topology,
Homotopy Types and the Fundamental Groupoid}, Ellis Horwood
Limited 1988, and of \textit{Elements of Modern Topology},
McGraw-Hill 1968.


\bibitem[BP00]{BuPa00} V.M.~Buchstaber, T.E.Panov, ~\textit{Torus
actions, combinatorial topology, and homological algebra}, Russian
Math.\ Surveys 55 (2000), no.\ 5, 825--921, arXiv:math.AT/0010073.

\bibitem[BP02]{BuPa02} V.M.~Buchstaber, T.E.Panov, ~\textit{Torus actions
and their applications in topology and combinatorics}, Univ.\
Lecture Ser., A.M.S.\ 2002.

\bibitem[CF60]{ConnerFloyd} P.E.~Conner, E.E.~Floyd, \textit{Fixed point free
involutions and equivariant maps}, Bull.\ Amer.\ Math.\ Soc.
\textbf{66} (1960), 416--441.



\bibitem[C95]{Connes} A.~Connes, \textit{Non commutative Geometry},
Academic Press 1995.


\bibitem[CrRy]{CrRy} H.~Crapo, J.~Ryan, \textit{Spacial realization of linear
scenes}, Structural Topology 13 (1986), 33--68.

\bibitem[DSS86]{DSS86} N.P.~Dolbilin, M.A.~Shtan'ko, and
M.I.Shtogrin, \textit{Cubic subcomplexes in regular lattices},
Dokl.\ Akad.\ Nauk., SSSR 291(1986), English translation:
\textit{Soviet Math.\ Dokl.}, 34(1987).

\bibitem[DSS87]{DSS87} N.P.~Dolbilin, M.A.~Shtan'ko, and M.I.Shtogrin.
\textit{Cubic manifolds in lattices}, Izv.\ Ross.\ Akad.\ Nauk.,
Ser.\ Mat.\ 58(1994), 93--107, English translation: \textit{Russ.\
Acad.\ Sci.\ Izv.\ Math.}, 44 (1995), 301--313.


\bibitem[DN02]{NovDyn} I.A.~Dynnikov, S.P.~Novikov, \textit{Geometry of the triangle
equation on two-manifolds}, arXiv:math-ph/0208041 v2 Oct 2002.


\bibitem[Epp99]{Epp99} D.~Eppstein, \textit{Linear complexity hexaedral mesh
generation}, Comput.\ Geom.\ 12 (1999), 3--16.

\bibitem[Fu99]{Fu99} L.~Funar, \textit{Cubulations, immersions, mappability and
a problem of Habegger}, Ann.\ Sci.\ E.N.S.\ 32 (1999), 681--700.

\bibitem[Fu99b]{Fu99b} L.~Funar, \textit{Cubulations mod bubble moves}, in
Proc.\ Conf.\ Low Dimensional Topology, Funchal, Madeira 1998
(H.~Nencka, Ed.) \textit{Contemporary Math.} 233, 29--43, A.M.S.\
1999.

\bibitem[Fu05]{Fu05} L.~Funar, \textit{Surface quadrangulations mod flips},
preprint (January, 2005) available at
\url{http://www-fourier.ujf-grenoble.fr/~funar}.

\bibitem[Ga04]{Gaif} A.A.~Gaifulin, \textit{Local formulae for combinatorial
Pontriagin classes}, Izvestya RAN: Ser.\ Mat.\ 68:5 13--66.


\bibitem[GKM]{GoKottMcPh} M. Goresky, R. Kottwitz, R. MacPherson,
\textit{Equivariant cohomology, Koszul duality, and the
localization theorem}, Invent.\ Math.\ 131 (1998), 25--83.

\bibitem[H71]{Higg}
P.J.~Higgins, \textit{Categories and groupoids}, Van Nostrand
Reinhold Co., London 1971.


\bibitem[I01]{Izm2001} I.~Izmestiev, \textit{$3$-dimensional manifolds defined
by simple polytopes with colored facets}, Mat.\ Zametki
\textbf{69} (2001), 375--382,  (English translation: Math.\ Notes
\textbf{69}, 1987).

\bibitem[I01b]{Izm-Usp2001} I.~Izmestiev, \textit{Free actions of a torus on a
manifold $\mathcal{L}_P$ and the projectivity group of a polytope
$P$}, Uspekhi Mat.\ Nauk, \textbf{56} (2001) pp.\ 169--170,
English translation: Russian Math.\ Survays \textbf{56} (2001).

\bibitem[IJ02]{IzmJos2002} I.~Izmestiev and M.~Joswig, \textit{Branched coverings,
triangulations, and $3$-manifolds}, Adv.\ Geom. \textbf{3} (2003),
191--225, arXiv:math.GT/0108202 v2 20 Mar 2002.

\bibitem[J01]{Josw2001} M.~Joswig, \textit{Projectives in simplicial
complexes and coloring of simple polytopes}, Math.\ Z.
\textbf{240} (2002) no.\ 2, 243--259, arXiv:math.CO/0102186 v3 27
Jun 2001.

\bibitem[J01b]{Jos-Usp} M.~Joswig, \textit{Projectives in simplicial
complexes and coloring of simple polytopes} (extended abstract),
Russ.\ Math.\ Surv.\ \textbf{56} (2001), 584--585.

\bibitem[K91]{Karalashvili} O.R.~Karalashvili, \textit{On mappings of
cubic manifolds into the standard lattice of Euclidean space},
Trudy Mat.\ Inst.\ Steklov 196(1991), 86--89, translated in:
Proc.\ Steklov Inst.\ Math.\ 196(1992).


\bibitem[K95]{Kirby} R.~Kirby, \textit{Problems in low-dimensional
topology}, in ``Geometric Topology'', Georgia International
Topology Conference, (W.H.~Kazez, Editor), AMS-IP Studies in
Advanced Math., $\mathbf{2}$, part $2$, 35--473, 1995.


\bibitem[Ko99]{Kozlov99} D.N.~Kozlov, \textit{Complexes of directed trees},
J.\ Combin.\ Theory Ser.\ A \textbf{88} (1999), no.\ 1, 112--122.


\bibitem[Ko04]{Kozl2004} D.N.~Kozlov, \textit{A simple proof for
folds on both sides in complexes of graph homomorphisms},
arXiv:math.CO/0408262 v2 Dec 2004.


\bibitem[Ko]{Kozlov-review} D.N.~Kozlov,
\textit{Chromatic numbers, morphism complexes, and Stiefel-Whitney
cha\-rac\-teristic classes}, "Geometric Combinatorics", IAS/Park
City Mathematics Series 14, American Mathematical Society,
Providence, RI; Institute for Advanced Study (IAS), Princeton, NJ.


\bibitem[Ko06]{Kozlov-short} D.N.~Kozlov,
\textit{Cobounding odd cycle colorings}, arXiv:math.AT/0602561 Feb
2006.

\bibitem[L78]{Lovasz} L. Lov\'{a}sz,
\textit{Kneser's conjecture, chromatic number and homotopy}, J.\
Comb.\ Theory, Ser. A, \textbf{25}:319--324, 1978.

\bibitem[M03]{Matousek} J. Matou\v{s}ek,
\textit{Using the Borsuk-Ulam Theorem; Lectures on Topological
Methods in Combinatorics and Geometry}, Springer Universitext,
Berlin 2003.

\bibitem[N96]{Novikov} S.P.~Novikov, \textit{Topology I},
Encyclopaedia Math.\ Sci.\ \textbf{12}, Springer-Verlag, Berlin
1996.

\bibitem[N04]{Novikov04} S.P.~Novikov, \textit{Discrete connections
on the triangulated manifolds and difference linear equations},
arXiv:math-ph/0303035 v2 Apr 2004.


\bibitem[Pe]{Pe} R.~Penrose, \textit{On the cohomology of impossible
figures}, Structural Topology 17 (1991), 11--16.

\bibitem[SZ04]{SZ04} A.~Schwartz, G.M.~Ziegler, \textit{Construction
techniques for cubical complexes, odd cubical $4$-polytopes, and
prescribed dual manifolds}, Experimental Math., 13:385–413, 2004,
arXiv:math.CO/0310269 v3 2 Jan 2004.


\bibitem[S05]{Schultz} C.~Schultz, \textit{A short proof of
$w_1^n(Hom(C_{2r+1},K_{n+2}))=0$ for all $n$ and a graph colouring
theorem of Babson and Kozlov}, arXiv:math.AT/0507346 v2 Aug 2005.

\bibitem[S06]{Schultz06} C.~Schultz, \textit{Graph colorings, spaces of edges
and spaces of circuits}, preprint 2006.

\bibitem[W96]{Wein96} A.~Weinstein, \textit{Groupoids: Unifying
Internal and External Symmetry}, Notices of the A.M.S., vol.\
\textbf{43}, July 1996.


\bibitem[Z92]{Zieg-text} G.M.~Ziegler, \textit{Some ``embedding
problems''}, unpublished manuscript December 1992.


\bibitem[Z98]{Zieg-book} G.M.~Ziegler, {\it Lectures on Polytopes},
Graduate Texts in Mathematics 152, Springer 1995, 2nd ed.\ 1998.


\bibitem[\v Z04]{Ziv04}
R. \v Zivaljevi\' c, \textit{Topological methods}, in CRC Handbook
of Discrete and Computational Geometry (new edition), J.E.
Goodman, J. O'Rourke (eds.), Boca Raton 2004.


\bibitem[\v Z05]{Ziv-groupoids} R. \v Zivaljevi\'c, \textit{Combinatorial
groupoids, cubical complexes, and the Lov\' asz conjecture},
arXiv:math.CO/0510204 v2 Oct 2005.

\bibitem[\v Z05b]{Ziv2005} R. \v Zivaljevi\'c, \textit{Parallel transport
of $Hom$-complexes and the Lov\' asz conjecture},
arXiv:math.CO/0506075 v1 3 Jun 2005.
\end{thebibliography}
\end{document}